\documentclass{article}
\usepackage{geometry}
\usepackage{authblk}
\usepackage{amsfonts}
\usepackage{bbold}
\usepackage{amssymb}
\usepackage{natbib}

\title{
Eigenfunctions of the Multidimensional\\
Linear Noise Fokker-Planck Operator via Ladder Operators
}

\author[1]{Todd K. Leen\footnote{Todd.Leen@georgetown.edu}}
\affil[1]{Graduate School of Arts and Sciences, Georgetown University}
\author[2]{Robert Friel}
\affil[2]{Courant Institute of Mathematical Sciences, New York, NY}
\author[3]{David Nielsen}

\begin{document}

\maketitle
  
\begin{abstract}
The eigenfunctions and eigenvalues of the Fokker-Planck operator with
linear drift and constant diffusion are required for expanding
time-dependent solutions and for evaluating our recent perturbation
expansion for probability densities governed by a nonlinear master
equation.  Although well-known in one dimension, for multiple
dimensions the eigenfunctions are not explicitly given in the
literature.  We develop raising and lowering operators for the
Fokker-Planck (FP) operator and its adjoint, and use them to obtain
expressions for the corresponding eigenvalues and eigenfunctions.  We
show that the eigenfunctions for the forward and adjoint FP operators
form a bi-orthogonal set, and that the eigenfunctions reduce to sums
of products of Hermite functions in a particular coordinate system.
\end{abstract}


\section{Introduction}
The Fokker-Planck equation (FPE) with linear drift and constant
diffusion describes an Ornstein-Uhlenbeck process.  In one dimension,
the eigenfunctions are the well-known Hermite functions
\citep{Risken:FokkerPlanck,Gardiner:Handbook}.  The eigenfunctions
enable expansion of time-dependent solutions to the FPE, and are
required to evaluate our recent perturbation expansion for the
probability densities arising from a {\em nonlinear} master equation
\citep{leenFriel2011,leenFriel2012,leenFrielNielsen2012}.
(Thomas and Grima have recently
derived a similar approximation and applied it to one-dimensional
chemical and gene expression systems \citeyearpar{ThomasGrima2015}.)  The
multi-dimensional eigenfunctions are {\em not} given in convenient form
in the literature.  (Liberzon and Brockett
\citeyearpar{LiberzonBrockett:SIAM2000} discuss the eigenvalue
spectrum, but do not give an explicit form of the complete set of
eigenfunctions, nor discuss the eigenfunctions of the adjunct operator.)
In the restricted case that the drift and diffusion
are simultaneously diagonalizable, the equations separate and the
solution reduces to products of the one-dimensional eigenfunctions.
However for general drift and diffusion, the eigenfunctions are not
simply obtained from the differential operators.

This note gives raising and lowering operators for the forward and
adjoint Fokker-Planck operators that develop the corresponding
complete sets of eigenfunctions.  We give bi-orthogonality relations
and show that in a special coordinate system, the eigenfunctions are sums
of products of Hermite functions.

\section{Operators and Eigenfunctions}
Let $x \in \mathbb{R}^N$, $\partial_i\equiv \partial/\partial x^i$,
$i=1,\ldots,N$.  
The forward and adjoint (backward) Fokker-Planck
operators for an N-dimensional Ornstein-Uhlenbeck process are
\begin{equation}
L \, q(x) \;\equiv\; -\, \partial_i \left(A^i_{\;j}\, x^j \, q \right)
                      \;+\; \frac{1}{2} B^{ij} \, \partial_i \, \partial_j \,q
\label{L}
\end{equation}
and
\begin{equation}
L^\dagger \, q(x) \;\equiv \; A^i_{\;j} \, x^j \; \partial_i \, q
                      \;+\; \frac{1}{2} B^{ij} \, \partial_i \, \partial_j \,q
\label{adjointL}
\end{equation}
where repeated upper and lower indices in {\em lower case Latin} are
summed over.  Thus \linebreak $(Ax)^i = A^i_{\;j}\,x^j \equiv \sum_j
A^i_{\,j} \, x^j$ is the $i^{th}$ component of the (linear) drift
vector.  We assume the diffusion matrix, whose elements are $B^{ij}$,
is positive definite.

Throughout we assume that the left $w_I$ and right $e_I$ eigenvectors of
$A$ form a complete set for $\mathbb{R}^n$, and that their corresponding
eigenvalues (perhaps complex) have negative real part\footnote{Since
$A$ is real, any complex eigenvectors occur as conjugate pairs.}.
These eigenvectors are normalized and bi-orthogonal
\begin{equation}
w_I \cdot e_J \;=\; w^*_I \cdot e^*_J \;=\; \delta_{IJ} \;\;.
\label{ewOrtho}
\end{equation}

\subsection{Stationary States}

The stationary state satisfies
\begin{displaymath}
L \; f_0(x) \;=\; 0 \;\;.
\end{displaymath}
and is given by
\begin{equation}
f_0(x) \;=\; \frac{1}{\sqrt{(2 \pi)^N \mbox{det}\Sigma}} 
  \, \exp \left(-\frac{1}{2} \, x^T \Sigma^{-1} x\right) \;\;\;.
\label{f0}
\end{equation}
The covariance matrix $\Sigma$ is a solution to the Liapunov equation
\begin{equation}
A\, \Sigma \;+\; \Sigma \, A^T \;=\; -B \;\;\;.
\label{Liapunov}
\end{equation}
The corresponding stationary state of $L^{\dagger}$ satisfies
\begin{displaymath}
L^{\dagger} \, g_0 \;=\; 0 
\end{displaymath}
and is given by
\begin{equation}
g_0(x) \;=\; 1 \;\;\;.
\end{equation}
As we will show, the eigenfunctions of $L$ and $L^\dagger$ form a
biorthogonal set, and for the stationary (ground) states we have
\begin{displaymath}
\left<g_0,f_0 \right> \;\equiv\; \int g_0^*(x) \, f_0(x) \; d^Nx \;=\; 1 
\end{displaymath}
as follows from the normalization of the multidimensional Gaussian $f_0$.

\subsection{Raising Operators}
We generate the complete set of eigenfunctions of $L$ and 
$L^{\dagger}$ by application of raising operators to $f_0$ and $g_0$.  
The next two subsections address the two sets of eigenfunctions in turn.

\subsubsection{Forward Eigenfunctions}
Let $e_I$, $I=1,\ldots,N$ be the right eigenvectors of $A$
with corresponding eigenvalues $\lambda_I$.  (Recall we assume that all
the eigenvalues of $A$ have negative real part, corresponding to an
asympotitically-stable fixed point at $x=0$.)  Let $e_I^i$ denote the $i^{th}$ 
component of the $I^{th}$ eigenvector.  The operators
\begin{equation}
V_I \, \;=\; -e_I^i \, \partial_i \;\equiv\; -e_I \cdot \nabla \;\;, \;\;\;
   I=1,2,\ldots N
\label{fRaise}
\end{equation}
satisfy the commutation relation
\begin{equation}
\left[L,\, V_I  \right] \;=\; \lambda_I \, V_I \;\;.
\label{fRaiseComm}
\end{equation}
(There is no sum on $I$.)
The commutator (\ref{fRaiseComm}) establishes $V_I$ as a raising
operator for eigenfunctions of $L$: If $f_{\lambda}$ is an
eigenfunction of $L$ with eigenvalue $\lambda$, then $V_I f_\lambda$
is also an eigenfunction of $L$ with eigenvalue $\lambda+\lambda_I$.
The application of products of the $V_I$ to $f_0$ generate new
eigenfunctions of $L$ which we denote by subscripts indicating the
number of applications of each operator. For example
\begin{displaymath}
f_{n_1,n_2,0,0 \ldots} \; \equiv \; (V_1)^{n_1} \, ( V_2)^{n_2} \, f_0 
\;\; \mbox{  with  } \;\;
 L \, f_{n_1,n_2,0,0 \ldots} \;=\; 
             (n_1 \lambda_1 + n_2 \lambda_2) \, f_{n_1,n_2,0,0 \ldots}
\end{displaymath}
Since $[V_I,V_J]=0$, the subscript labels on $f$ uniquely determine the 
eigenfunctions.  We build a general eigenfunction of $L$ as
\begin{equation}
f_{n_1,n_2,\ldots,n_N} \; \equiv \; \left(\prod_{I=1}^N \, (V_I)^{n_I} \right) \; f_0 
\label{generalEigF}
\end{equation}
which satisfies
\begin{equation}
L \; f_{n_1,n_2,\ldots,n_N} \;=\;  \left(\sum_{I=1}^N n_I\lambda_I \right) 
          \, f_{n_1,n_2,\ldots,n_N} \;\;.
\label{generalValF}
\end{equation}
Clearly from Eqn.\ (\ref{generalEigF})
\begin{equation}
(V_I)^k \, f_{n_1,\ldots,n_I,\ldots,n_N} \;=\; f_{n_1,\ldots,(n_I+k),\ldots,n_N} \;\;.
\end{equation}

\subsubsection{Adjoint Eigenfunctions}
Let $w_I$, $I=1,\ldots,N$ be the left eigenvectors of $A$ with
eigenvalues $\lambda_I$, and let $w_{I\,i}$ denote its $i^{th}$ component.
(The {\em subscript} indicating the component is appropriate since we
regard $w_I$ as a {\em co-vector}.) Define the operators
\begin{equation} \bar{V}_I \;\equiv\; 2(w^*_{I\,i} \, x^i) \;+\;
 2 \left[(A+\lambda_I)^{-1} B w_I^* \right]^i \partial_i \;=\; 
 2\left[ \, w^*_{I\,i} \, x^i \;-\; (\Sigma w_I^*)^i \partial_i \, \right]\;\;,
\label{bRaise} 
\end{equation} where the second equality
follows from the Liapunov equation (\ref{Liapunov}).
This operator satisfies the commutation relation
\begin{equation}
\left[L^\dagger, \bar{V}_I\right] \;=\; \lambda_I^* \bar{V}_I \;,
\label{bRaiseOpComm}
\end{equation}
which establishes it as a raising operator for eigenfunctions
of $L^\dagger$:  If
$g_\lambda$ is an eigenfunction of $L^\dagger$ with 
eigenvalue $\lambda$, then $\bar{V}_I g_\lambda$ is an eigenfunction with
eigenvalue $\lambda+\lambda_I^*$.
Analogously to the forward eigenfunctions, acting with $\bar{V}_I$ on $g_0$
generates a new eigenfunction which we denote by subscripts
indicating the number of applications of each raising operator.
This subscript notation
is free of ambiguity since the raising operators
corresponding to different eigenvectors of $A$ commute
\begin{displaymath}
\left[\bar{V}_I,\bar{V}_J \right] \;=\; -2 (w^*_J)^T \Sigma w_I^* 
                                      + 2 (w_I^*)^T \Sigma w^*_J \;=\; 0 \;\;.
\end{displaymath}
We construct a general eigenfunction of $L^\dagger$ by repeated 
application of raising operators.  Thus
\begin{equation}
g_{n_1,n_2,\ldots,n_N} \;=\; \prod_{I=1}^N \left(\bar{V}_I\right)^{n_I} \; g_0 \;,
\label{generalEigB}
\end{equation}
which satisfies
\begin{equation}
L^\dagger \; g_{n_1,n_2,\ldots,n_N} \;=\; \left(\sum_{I=1}^N n_I \lambda_I^* \right) 
                                              g_{n_1,n_2,\ldots,n_N} \;\;.
\label{generalValB}
\end{equation}
Clearly, from Eqn.\ (\ref{generalEigB})
\begin{equation}
\left(\bar{V}_I\right)^k \, g_{n_1,\ldots,n_I,\ldots,n_n} \;=\;
   g_{n_1,\ldots,(n_I+k),\ldots,n_n} \;\;.
\end{equation}

\subsection{Lowering Operators}
The two sets of raising operators are complemented by lowering operators.  Their
effect on the states, derived here, makes proving the bi-orthogonality
property trivial.  

\subsubsection{Adjoint Eigenfunctions} 
Taking the adjoint of  Eqn.\ (\ref{fRaiseComm}) yields
\begin{equation}
\left[ L^\dagger,V_I^\dagger \right] \;=\; -\lambda_I^* \, V_I^\dagger
\label{lowerG}
\end{equation}
which establishes $V_I^\dagger = e_I^* \cdot \nabla$ 
as a lowering operator for eigenfunctions
of $L^\dagger$: If $g_\lambda$ is an eigenfunction of $L^\dagger$ with eigenvalue
$\lambda_I$, then $V_I^\dagger \, g_\lambda$ is an eigenfunction of $L^\dagger$
with eigenvalue $\lambda - \lambda_I^*$.  
In particular, $V_I^\dagger$ \em{kills the stationary state}
\begin{displaymath}
V_I^\dagger \, g_0 = e_I^* \cdot \nabla \, g_0 =0 \;\;.
\end{displaymath}
Using the commutators
\begin{equation}
\left[ V_I^\dagger ,\bar{V}_J \right] \;=\; 2\, e_I^* \cdot w_J^* \;=\;
           2\, \delta_{IJ}
\label{commuteVVbar}
\end{equation}
and $[\bar{V}_I,\bar{V}_J]=0$, and the expression for the adjoint eigenfunctions
Eqn.\ (\ref{generalEigB}),
the action of $V_I^\dagger$ on eigenfunctions of $L^\dagger$ follows as
\begin{eqnarray}
V_J^\dagger \; g_{n_1,n_2,\ldots,n_N} &\equiv& V_J^\dagger \, \left(  \prod_{I=1}^N \, 
          \left( \bar{V}_I\right)^{n_I} \right) \, g_0 
\;=\;  \left( \prod_{I\ne J} \, \left( \bar{V}_I\right)^{n_I}\right)
             V_J^\dagger \, \left(\bar{V}_J \right)^{n_J} \, g_0 \nonumber \\
&=& \left( \prod_{I \ne J} \, \left( \bar{V}_I \right)^{n_I} \right)
  \left[\left(\bar{V}_J\right)^{n_J} \, V_J^\dagger 
  \,+\, (2n_J) \left(\bar{V}_J\right)^{n_J-1} \right] \, g_0 \nonumber \\[1.25pc]
&=& (2 n_J) \, g_{n_1,n_2,\ldots,n_J-1,\ldots,n_N} 
\label{gLower}
\end{eqnarray}
where the second line follows from the first by commuting
$V_J^\dagger$ past all $n_J$ factors of $\bar{V}_J$, and the third
line follows since $V_J^\dagger \, g_0\equiv e_J^* \cdot \nabla \, g_0 = 0$.
Multiple applications yield
\begin{equation}
\left(V_J^\dagger \right)^k \; g_{n_1,n_2,\ldots,n_N} 
   \;=\; \left\{ 
        \begin{array}{cc}
        2^k\frac{n_J!}{(n_J-k)!}  g_{n_1,n_2,\ldots,n_J-k,\ldots, n_N} \;, & k \le n_J \\
            0\;, & k > n_J 
        \end{array}\right. \;\;.
\label{gLowerMany}
\end{equation}

\subsubsection{Forward Eigenfunctions}
Taking the adjoint of Eqn.\ (\ref{bRaiseOpComm}) we recover
\begin{equation}
\left[\, L, \bar{V}_I^\dagger \, \right] \;=\; 
     -\lambda_I \, \bar{V}_I^\dagger
\label{lowerF}
\end{equation}
which establishes $\bar{V}_I^\dagger$ as a 
lowering operator for eigenfunctions of $L$: If $f_\lambda$ is an eigenfunction
of $L$ with eigenvalue $\lambda$, then $\bar{V}_I^\dagger \, f_\lambda$ is an
eigenfunction with eigenvalue $\lambda-\lambda_I$.  In particular, it 
\em{kills the stationary state}
\begin{equation}
\bar{V}_I^\dagger \, f_0 \;=\; 2\left( w_{I\,i} \, x^i + (\Sigma w_I)^i 
             \partial_i \right) \, f_0 \;=\; 0
\label{killf0}
\end{equation}
having used the definition of $f_0$ in Eqn.\ (\ref{f0}).  
Similarly to the
case for the adjoint lowering operators, it is straightforward to show that
\begin{equation}
\bar{V}_J^\dagger \, f_{n_1,n_2,\ldots,n_N} \;=\; (2n_J) \, 
    f_{n_1,\ldots,n_J-1,\ldots,n_N}  \;\;.
\label{fLowerGen}
\end{equation}
Multiple applications yield
\begin{equation}
(\bar{V}_J^\dagger)^k \, f_{n_1,\ldots,n_N} \;=\; 
  \left\{ \begin{array}{cc}
      2^k \frac{n_J!}{(n_J-k)!} \, f_{n_1,\ldots,n_J-k,\ldots,n_N} \;, & k \le n_J \\
                   0 \;, & k > n_J \
          \end{array} \right. \;\;.
\label{fLowerMany}
\end{equation}

\section{Bi-Orthogonality and Normalization}

The two sets of functions form a bi-orthogonal set.  The usual result,
that eigenfunctions corresponding to different eigenvalues are
orthogonal can be strengthened.  
We will show that
\begin{eqnarray}
\left< \, g_{m_1,\ldots m_N}, \; f_{n_1, \ldots, n_N} \, \right>
&\equiv& \int \, g^*_{m_1,\ldots m_N}(x) \; f_{n_1, \ldots, n_N}(x) \;\; d^Nx \nonumber \\
&=&
\left( \prod_{I=1}^N \delta_{m_I,n_I} \;2^{n_I}\,(n_I!) \right)  \; <g_0,f_0>
\nonumber \\
&=&
\left( \prod_{I=1}^N \delta_{m_I,n_I} \;2^{n_I}\,(n_I!) \right)
   \;\;.
\label{orthonormal}
\end{eqnarray}
We consider three cases:

\paragraph{Case I} Suppose $n_J > m_J$ for some $J$.  Write
\begin{eqnarray}
\left<g_{m_1,\ldots,m_J,\ldots,m_N},f_{n_1,\ldots,n_J,\ldots,n_N}\right> &=&
\left<g_{m_1,\ldots,m_J,\ldots,m_N},\, V_J^{n_J} \, f_{n_1,\ldots,n_J=0,\ldots,n_N}\right> 
\nonumber \\
&=& \left<\, (V_J^\dagger)^{n_J} \, g_{m_1,\ldots,m_J,\ldots,m_N} ,
                           \, f_{n_1,\ldots,n_J=0,\ldots,n_N} \right> \;=\; 0
\label{caseI}
\end{eqnarray}
having used Eqn.\ (\ref{gLowerMany}) to arrive at the last equality.

\paragraph{Case II}
Now suppose instead that $m_J > n_J$ for some $J$.  Write
\begin{eqnarray}
\left<g_{m_1,\ldots,m_J,\ldots,m_N},f_{n_1,\ldots,n_J,\ldots,n_N}\right> &=&
\left<\, \bar{V}_J^{m_J}\, g_{m_1,\ldots,m_J=0,\ldots,m_N},
      \, f_{n_1,\ldots,n_J,\ldots,n_N} \right> \nonumber \\
&=& \left<\, g_{m_1,\ldots,m_J=0,\ldots,m_N} ,
       \, (\bar{V}_J^\dagger )^{m_J} f_{n_1,\ldots,n_J,\ldots,n_N} \right> \;=\; 0
\label{caseII}
\end{eqnarray}
having used Eqn.\ (\ref{fLowerMany}) to arrive at the last equality.

\paragraph{Case III} Suppose $n_J = m_J$ for all $J$.  Write
\begin{eqnarray}
\left<g_{n_1,\ldots,n_N},f_{n_1,\ldots,n_N}\right> &=&
\left<g_{n_1,\ldots,n_N},\, V_1^{n_J} \, f_{0,n_2,\ldots,n_N}\right> \nonumber \\
&=& \left<\, (V_1^\dagger)^{n_1} \, g_{n_1,n_2,\ldots,m_N} ,
                           \, f_{0,n_2,\ldots,n_N} \right> \nonumber \\
&=&  2^{n_1}\,n_1! \; \left<g_{0,n_2,\ldots,n_N} ,
                           \, f_{0,n_2,\ldots,n_N} \right> \nonumber \\
&=&  2^{n_1}\,n_1! \; \left< (V_2^\dagger)^{n_2} \, g_{0,n_2,n_3,\ldots,n_N} ,
                           \, f_{0,0,n_3,\ldots,n_N} \right> \nonumber \\
&=& 2^{n_1}\, 2^{n_2}\,(n_1!)(n_2!) \; \left< g_{0,0,n_3,\ldots,n_N} ,
                           \, f_{0,0,n_3,\ldots,n_N} \right> \nonumber \\
& \vdots & \nonumber \\
&=& \left< \prod_{J=1}^N \, 2^{n_J} (n_J!) \right> \, \left< g_0,f_0 \right>
                           \nonumber \\
&=& \left< \prod_{J=1}^N \, 2^{n_J} (n_J!) \right>
\label{caseIII}
\end{eqnarray}
having used Eqn.\ (\ref{gLowerMany}) repeatedly.  The three cases together prove
the desired result Eqn.\ (\ref{orthonormal}).

\section{Analytic Form of the Eigenfunctions}
In general coordinates, the eigenfunctions
$f_{n_1,n_2,\ldots,n_N}$ and $g_{n_1,\ldots,n_N}$ do not assume a familiar
analytic form.  However they do if we transform to coordinates in
which $\Sigma = \frac{1}{2} \, \mathbb{1}$.  (Since $\Sigma$ is a
real, symmetric, positive-definite matrix, this is always possible.)
In these coordinates,
\begin{eqnarray}
f_0(x) &=&  \frac{1}{\sqrt{\pi^N}} \, e^{-x^T\,x} 
\label{f0special} \\
g_0(x) &=& 1 \;\;.
\label{g0special}
\end{eqnarray}
In these coordinates the forward eigenfunctions
(\ref{generalEigF}) are sums of products of Hermite polynomials
times $f_0$ and the backward eigenfunctions (\ref{generalEigB})
are sums of products of Hermite polynomials.

\subsection{Forward Eigenfunctions}
To start, note that the usual generating expression for the Hermite
polynomials \citep{AbramowitzStegan} can be rearranged to read
\begin{equation}
H_n(x_i) \, e^{-x^Tx} \;=\; (-\partial_i)^n e^{-x^Tx}
 \;\;\;.
\label{HermiteGen}
\end{equation}
Hence, application of $V_I$ to $f_0$ yields
\begin{displaymath}
V_I \, f_0(x) \;=\; -e_I^i \partial_i \; f_0(x) \;=\; 
 \sum_{i=1}^N e_I^i \, H_1(x_i) \, f_0(x) 
\end{displaymath}
a linear combination of Hermite functions in each of the variables
$x_1,\ldots,x_N$.  (When a coordinate index falls inside a function
argument --- as in the last expression --- we will write the summation
explicitly to avoid confusion.)
Applying $V_I$ to $f_0$, $n_I$ times results in
\begin{displaymath}
V_I^{n_I} f_0 \;=\;  \left( -e^i_I \partial_i \right)^{n_I} \, f_0
\end{displaymath}
which can be evaluated using the multinomial theorem.  Explicitly
\begin{eqnarray}
V_I^{n_I} f_0 &=& \sum_{(i_1+\ldots+i_N=n_I)}
           \left(\frac{n_I!}{i_1! \, i_2! \cdots i_n!}\right)
           \prod_{1\le k \le N} (-e^k_I \, \partial_k)^{i_k} \; f_0 \nonumber \\
           &=& \sum_{(i_1+\ldots+i_N=i_I)}
           \left(\frac{n_I!}{i_1! \, i_2! \cdots i_n!}\right)
           \prod_{1\le k \le N} (e^k_I)^{i_k} H_{i_k}(x_k) \; f_0(x)  \nonumber
\label{forwardVn2}
\end{eqnarray}
where the summation is over all values of indices satisfying the
constraint $i_1+i_2+\cdots+i_N=n_I$.

The action of two distinct ladder operators multiple times is
\begin{eqnarray}
V_I^{n_I} \, V_J^{n_J} f_0 &=& 
\sum_{\footnotesize{\begin{array}{c}(i_1+\ldots+i_N=n_I)\\(j_1+\ldots+j_N=n_J)\end{array}}}
\left( \frac{n_I!}{i_1! i_2! \ldots i_N!} \right)
\left( \frac{n_J!}{j_1! j_2! \ldots j_N!} \right) 
  \prod_{1\le k \le N}(e_I^k)^{i_p} \, (e_J^k)^{j_k} (-\partial_k)^{i_k+j_k} \, f_0
  \nonumber \\[0.5pc]
&=&
\sum_{\footnotesize{\begin{array}{c}(i_1+\ldots+i_N=n_I)\\(j_1+\ldots+j_N=n_J)\end{array}}}
\left( \frac{n_I!}{i_1! i_2! \ldots i_N!} \right)
\left( \frac{n_J!}{j_1! j_2! \ldots j_N!} \right) 
  \prod_{1\le k \le N}(e_I^k)^{i_p} \, (e_J^k)^{j_k} H_{i_k+j_k}(x_k) \, f_0
 \nonumber \\
\end{eqnarray}
\normalsize
This generalizes in the obvious way to products of the form
\begin{equation}
f_{n_1,n_2,\ldots,n_N}(x) \;=\; V_1^{n_1} \, V_2^{n_2} \ldots V_N^{n_N} \, f_0(x) \;\;.
\end{equation}
So in our special coordinates the eigenfunctions of $L$ are 
{\em sums of products of Hermite polynomials} times $f_0$.

\subsection{Adjoint Eigenfunctions}
In our special coordinates $\Sigma^{ij} = \frac{1}{2} \delta^{ij}$
so the backward raising operator 
(\ref{bRaise}) simplifies to
\begin{equation}
\bar{V}_I \;=\; 2\, (w^*_{I\,i} \, x^i) \,-\, (w_I^*)^i \, \partial_i \;\;,
\label{bLadderOpSpecial}
\end{equation}
where $(w_I^*)^i \equiv \delta^{ij}w^*_{I\,j}$.  
From the standard recursion relations for the Hermite polynomials one
has
\begin{equation}
\left(2x-\frac{d}{dx} \right) \, H_n(x) \;=\; H_{n+1}(x) \;\;.
\end{equation}
Thus using the form of $\bar{V}_I$ in our special coordinates 
(\ref{bLadderOpSpecial}), its action on a product of Hermite
polynomials is
\begin{displaymath}
\bar{V}_I \; H_{i_1}(x_1)\,H_{i_2}(x_2)\,\cdots\,H_{i_N}(x_N) \;=\;
          \sum_{k=1}^N \, 
            H_{i_1}(x_1)\,\,\cdots\,w^*_{I\,k} \, H_{i_k+1}(x_k) \,\cdots \,
                        H_{i_N}(x_N) \;\;.
\end{displaymath}
It is convenient to define the operator $r^k$ which raises the order
of $H_{m}(x_k)$ by unity,
\begin{displaymath}
r^k \; H_{i_1}(x_1) \ldots H_{i_k}(x_k) \ldots H_{i_N}(x_N) 
     \;=\;  H_{i_1}(x_1) \ldots H_{i_k+1}(x_k) \ldots H_{i_N}(x_N) \;\;\;.
\end{displaymath}
Then the action of $\bar{V}_I$ on a product of Hermite polynomials is
\begin{eqnarray}
\bar{V}_I \; H_{i_1}(x_1)\ldots H_{i_N}(x_N) &=&
     \sum_{k=1}^N w^*_{I\,k} \;  r^k \; H_{i_1}(x_1)\ldots H_{i_N}(x_N) \nonumber \\
  &=& \sum_{k=1}^N  \,H_{i_1}(x_1)\ldots w^*_{I\,k} \, H_{i_k+1}(x_k) 
                    \ldots H_{i_N}(x_N)  \;.
\label{actionVBar2}
\end{eqnarray}
Writing
\begin{displaymath}
g_0 \;=\; 1 \;=\; H_0(x_1)\,H_0(x_2)\,\cdots\,H_0(x_N)
\end{displaymath}
and using Eqn.\ (\ref{actionVBar2}) gives
\begin{equation}
\bar{V}_I \, g_0 \;=\; \sum_{k=1}^N w^*_{I\,k} \, H_1(x_k) 
\label{barVonG0}
\end{equation}
The action of $\bar{V}_I^{n_I}$ on $g_0$ is evaluated
using the multinomial theorem
\begin{eqnarray}
\bar{V}_I^{n_I} \, g_0
 &=& \sum_{\footnotesize (i_1+\cdots+i_n = n_I)} 
     \left(\frac{n_I!}{i_1! \ldots i_N!}\right)
             \prod_{1 \le k \le N} (w^*_{I\,k} \; r^k)^{i_k} \; g_0
        \nonumber \\
 &=&  \sum_{\footnotesize (i_1+\cdots+i_n = n_I)} 
       \left(\frac{n_I!}{i_1! \ldots i_N!}\right)
           \prod_{1 \le k \le N} (w^*_{I\,k})^{i_k} H_{i_k}(x_k) \nonumber
\label{barVn}
\end{eqnarray}
The action of two such operators on $g_0$ is
\begin{eqnarray}
\bar{V}_J^{n_J} \; \bar{V}_I^{n_I} \, g_0
&=& \sum_{\footnotesize \begin{array}{c}(i_1+\cdots+i_N = n_I)\\
                                        (j_1+\cdots+j_N = n_J)
                        \end{array}} 
   \left(\frac{n_J!}{j_1! \ldots j_N!}\right)
   \left(\frac{n_I!}{i_1! \ldots i_N!}\right)
             \prod_{1\le k \le N}(w^*_{I\,k})^{i_k}(w^*_{J\,k})^{j_k} \; 
             H_{i_k+j_k}(x_k) \;.
\nonumber \\
\label{varVng0two}
\end{eqnarray}
This generalizes in the obvious way to evaluate the general backward 
eigenfunction in Eqn.\ (\ref{generalEigB}).

\section{Examples}
\label{examples}
\paragraph{Example 1 --- Derivative and Multiplication Operators}
Recall that the left ($w_I$) and right ($e_I$) 
eigenvectors of $A$ form a complete, bi-orthogonal set for $\mathbb{R}^N$
(see Eqn.\ (\ref{ewOrtho})).
Then, from the definition of the raising operator (\ref{fRaise}) for
eigenfunctions of $L$, we have
\begin{equation}
 \nabla \;=\; \sum_{I=1}^N \,  w_I^* \; V_I^\dagger \;\;.
\label{gradAsVs}
\end{equation}
Similarly, from the definition of the raising operator (\ref{bRaise}) for
eigenfunctions of $L^\dagger$ and this last result (\ref{gradAsVs}), 
we have
\begin{equation}
 x \;=\; \frac{1}{2} \sum_{I=1}^N e_I^*
        \left[ \bar{V}_I + 2 \sum_{J=1}^N
                (w_I^* \cdot \Sigma w^*_J) \, V_J^\dagger  \right] \;\;.
\end{equation}

\vspace{1pc}
\paragraph{Example 2 --- Forward and Adjoint FPE Operators} Using the results in Example 1, 
we can rewrite $L$ as
\begin{equation}
L \;=\; \frac{1}{2} \sum_I \, \lambda_I \, V_I \, 
                           \bar{V}_I^\dagger
\label{L0AsVs}
\end{equation}
which can be verified by its action on the eigenfunctions
\begin{eqnarray}
L \, f_{n_1,\ldots,n_N} &=&
 \frac{1}{2} \sum_I \, \lambda_I \, V_I \, \bar{V}_I^\dagger \;
   f_{n_1,\ldots,n_N}     
\;=\;   \sum_I \, 
         n_I \lambda_I \; V_I  f_{n_1,\ldots,n_I-1,\ldots,n_N} \nonumber \\
&=&  \sum_I \, n_I \lambda_I \; f_{n_1,\ldots,n_N} \;\;.
\end{eqnarray}
Similarly, we recover $L^\dagger$ as
\begin{equation}
L^\dagger \;=\; \frac{1}{2} \sum_I \lambda_I^* \, \bar{V}_I \, V_I^\dagger
\label{L0DaggerAsVs}
\end{equation}
which can be verified by its action on the eigenfunctions
$g_{n_1,\ldots,n_N}$.  (Both Eqns. (\ref{L0AsVs}) and (\ref{L0DaggerAsVs}) 
have the flavor of quantum oscillator Hamiltonians involving products
of lowering and raising operators\footnote{Hence, 
$\frac{1}{2} V_I \bar{V}_I^\dagger$ and $\frac{1}{2} \bar{V}_I \, V_I^\dagger$
are occupation number operators for the $I^{th}$
eigenstate in $f_K$ and $g_K$ respectively.}
The difference here is that $L$
and $L^\dagger$ are not Hermitian.  The correspondence with quantum
mechanics (in the 1-D case) is discussed by Gardiner
\citeyearpar{Gardiner:Handbook}.)

\vspace{1pc}
\paragraph{Example 3 --- Fourier Series} Time-dependent solutions
can be expanded in terms of the forward eigenfunctions exactly as
in the one-dimensional case.  Let $K\equiv\{k_1,k_2,\ldots,k_N\}$ 
denote an index set for the eigenfunctions and eigenvalues.  Then
\begin{equation}
F(x,t) \;=\; \sum_K \, \alpha_K \, f_K(x) \, \exp(\lambda_K \, t)
\label{expansion}
\end{equation}
clearly satisfies 
\begin{displaymath}
\partial_t \, F(x,t) \;=\; L \, F(x,t) 
\end{displaymath}
where the coefficients are given by $\alpha_K =\left<g_K,F(x,0)\,\right>$.
(Since the eigenvalues of $A$ have negative real part, all the exponentials
in (\ref{expansion}) are decaying.)
In multiple dimensions, the drift Jacobian $A^i_{\,j}$ can have
one or more pairs of complex-conjugate eigenvalues and
Eqn.\ (\ref{expansion}) can represent damped oscillating solutions.

\paragraph{Example 4 --- Solutions to Inhomogeneous Equations}  
In our perturbation solution for
densities satisfying a nonlinear master equation
\citep{leenFriel2011,leenFriel2012,leenFrielNielsen2012},
the $i^{th}$ order corrections to the 
equilibrium density is denoted $P^{(i)}(x)$ 
and is given by inhomogeneous equations of the form
\begin{equation}
L \, P^{(i)}(x) \;=\; q^{(i)}(x) \;\;,
\label{inhomogeneous}
\end{equation}
where the $q^{(i)}$ are known. Let $K\equiv\{k_1,k_2,\ldots,k_N\}$ 
denote an index set for the eigenfunctions and eigenvalues.  Next, expand
$P^{(i)}$ in a linear combination of the eigenfunctions $f_K(x)$, substitute
that into (\ref{inhomogeneous}), take the inner product with $g_J(x)$, and use
bi-orthogonality relations (\ref{orthonormal}) to obtain
\begin{equation}
  P^{(i)}(x) \;=\; \sum_{K\ne 0} \,
  \frac{\left<g_K,q^{(i)} \right>}{\lambda_K \, \left< g_K, f_K \right>} \,
  f_K(x)
\label{inhomogeneousSolution}
\end{equation}
where the summation {\em excludes} $f_0$.  

\section{Conclusion}
We have provided raising and lowering operators to develop the
eigenfunctions (and their corresponding eigenvalues) of the forward
and adjoint multidimensional Fokker-Planck operators for the
Ornstein-Uhlenbeck process.
The eigenfunctions form a basis for
expanding solutions to the time-dependent Fokker-Planck equation, and
for a perturbation expansion of the densities arising from a nonlinear
master equation
\citep{leenFriel2011,leenFriel2012,leenFrielNielsen2012,ThomasGrima2015}.  

We gave bi-orthogonality and normalization results.  We
showed that in coordinates for which the covariance of the stationary
state is spherically symmetric with variance one-half, the
eigenfunctions reduce to sums of products of Hermite polynomials times
$f_0$.  In applications to time-dependent solutions of the Fokker-Planck
equation and to inhomogeneous equations (see Example 3 and Example 4 in
Section \ref{examples}) one assumes the eigenfunctions $f_K$ form a complete
set on ${\cal L}_2$.  The proof of completeness is similar to that used to
show completeness of the one-dimensional Hermite functions.

\paragraph{Acknowledgments}  This work was supported by NSF under
grant IIS-0812687.  The authors thank Gerardo Lafferriere
and Crispin Gardiner for their comments.

\vspace{1pc}
\bibliography{ouladders}

\begin{thebibliography}{8}
\expandafter\ifx\csname natexlab\endcsname\relax\def\natexlab#1{#1}\fi
\expandafter\ifx\csname url\endcsname\relax
  \def\url#1{\texttt{#1}}\fi
\expandafter\ifx\csname urlprefix\endcsname\relax\def\urlprefix{URL }\fi
\providecommand{\eprint}[2][]{\url{#2}}
\providecommand{\bibinfo}[2]{#2}
\ifx\xfnm\relax \def\xfnm[#1]{\unskip,\space#1}\fi
\bibitem[{Abramowitz and Stegun(1972)}]{AbramowitzStegan}
\bibinfo{author}{Abramowitz, M.}, \bibinfo{author}{Stegun, L.},
  \bibinfo{year}{1972}.
\newblock \bibinfo{title}{Handbook of Mathematical Functions}.
\newblock \bibinfo{publisher}{U.S. Department of Commerce, National Bureau of
  Standards}.
\bibitem[{Gardiner(2009)}]{Gardiner:Handbook}
\bibinfo{author}{Gardiner, C.}, \bibinfo{year}{2009}.
\newblock \bibinfo{title}{Stochastic Methods, A Handbook for the Natural and
  Social Sciences, Fourth Edition}.
\newblock \bibinfo{publisher}{Springer-Verlag}, \bibinfo{address}{Berlin}.
\bibitem[{Leen and Friel(2011)}]{leenFriel2011}
\bibinfo{author}{Leen, T.K.}, \bibinfo{author}{Friel, R.},
  \bibinfo{year}{2011}.
\newblock \bibinfo{title}{Perturbation theory for stochastic learning
  dynamics}, in: \bibinfo{booktitle}{Proceedings of the IJCNN 2011},
  \bibinfo{publisher}{IEEE Press}, \bibinfo{address}{San Jose, CA}.
\bibitem[{Leen and Friel(2012)}]{leenFriel2012}
\bibinfo{author}{Leen, T.K.}, \bibinfo{author}{Friel, R.},
  \bibinfo{year}{2012}.
\newblock \bibinfo{title}{Stochastic perturbation methods for
  spike-timing-dependent plasticity}.
\newblock \bibinfo{journal}{Neural Computation.} \bibinfo{volume}{24},
  \bibinfo{pages}{1109--1146}.
\bibitem[{Leen et~al.(2012)Leen, Friel and Nielsen}]{leenFrielNielsen2012}
\bibinfo{author}{Leen, T.K.}, \bibinfo{author}{Friel, R.},
  \bibinfo{author}{Nielsen, D.}, \bibinfo{year}{2012}.
\newblock \bibinfo{title}{Approximating distributions in stochastic learning}.
\newblock \bibinfo{journal}{Neural Networks} \bibinfo{volume}{32},
  \bibinfo{pages}{219--228}.
\bibitem[{Liberzon and Brockett(2000)}]{LiberzonBrockett:SIAM2000}
\bibinfo{author}{Liberzon, D.}, \bibinfo{author}{Brockett, R.W.},
  \bibinfo{year}{2000}.
\newblock \bibinfo{title}{Spectral analysis of {Fokker-Planck} and related
  operators arising from linear stochastic differential equation}.
\newblock \bibinfo{journal}{SIAM J. Control Optim.} \bibinfo{volume}{38},
  \bibinfo{pages}{1453--1467}.
\bibitem[{Risken(1989)}]{Risken:FokkerPlanck}
\bibinfo{author}{Risken, H.}, \bibinfo{year}{1989}.
\newblock \bibinfo{title}{The Fokker-Planck Equation}.
\newblock \bibinfo{publisher}{Springer-Verlag}, \bibinfo{address}{Berlin}.
\bibitem[{Thomas and Grima(2015)}]{ThomasGrima2015}
\bibinfo{author}{Thomas, P.}, \bibinfo{author}{Grima, R.},
  \bibinfo{year}{2015}.
\newblock \bibinfo{title}{Approximate probability distributions of the master
  equation}.
\newblock \bibinfo{journal}{Physical Review} \bibinfo{volume}{E92}.
\newblock \bibinfo{note}{DOI: 10.1103/PhysRevE.92.012120}.

\end{thebibliography}
\bibliographystyle{model2-names}

\end{document}